\documentclass[a4paper,12pt]{article}
\usepackage{amsmath, amssymb, amsthm}
\usepackage{geometry}
\usepackage{graphicx}
\usepackage{xcolor}
\usepackage{hyperref}

\newtheorem{theorem}{Theorem}[section]

\theoremstyle{definition}
\newtheorem{definition}[theorem]{Definition}

\theoremstyle{remark}
\newtheorem{remark}[theorem]{Remark}
\DeclareMathOperator   {\tr}     {tr}

\title{\textbf{Almost Differentially Nondegenerate Nijenhuis Operators}}
\author{Dinmukhammed Akpan\\Email: dinmukhammed.akpan@uni-jena.de}
\date{}
\begin{document}
\maketitle

\section{Introduction}
The Nijenhuis geometry is a new section of modern differential geometry, which
studies smooth manifolds with a tensor field $L$ of type $(1,1)$, i.e.,
manifolds with endomorphism field for which the Nijenhuis torsion $N_L$ is
identically equal to zero. Fundamental results in this direction have been
obtained by A.\, V.~Bolsinov, V.\, S.~Matveev, and A.\, Yu.~Konyaev in the
papers \cite{BMK1, BMK2, BMK3, BMK4, BMK5}. The vanishing of the Nijenhuis
tensor for operator fields is a necessary (but not sufficient) condition for
their integrability (i.e., the possibility of reduction to a constant form by a
change of coordinates).

The paper \cite{BMTM} contains a problem number 5.13, which can be reformulated
as follows: \textit{describe all two-dimensional Nijenhuis operators for which
the differential of the trace is not equal to zero}. This problem was considered in the work \cite{Akpan}. Also, in \cite{Akpan}, a theorem was proved on the
form of a two-dimensional Nijenhuis operator whose trace is the first
coordinate and the determinant has a Morse singularity with respect to the
second variable. In this paper, we obtain a generalization of this result to an
arbitrary dimension (Theorem~\ref{th2}).

\section{Main Part}

\begin{definition}
Let $M^n$ be a smooth $n$-dimensional manifold and let $L$ be a tensor field of
type $(1,1)$. Then the Nijenhuis torsion or the Nijenhuis tensor $N_L$ is a
tensor of type $(1,2)$ which is invariantly defined as follows:
$$
N_L[u,v] = L^2[u,v] + [Lu, Lv] - L[u, Lv] - L[Lu, v],
$$
where $u,v$ are arbitrary vector fields and $[u,v]$ is their commutator. In 
coordinates, this condition is written as follows:
$$
(N_L)^i_{jk} = L^l_j \frac{\partial L^i_k}{\partial x^l} - L^l_k \frac{\partial
L^i_j}{\partial x^l} - L^i_l \frac{\partial L^l_k}{\partial x^j} + L^i_l
\frac{\partial L^l_j}{\partial x^k}.
$$
\end{definition}

\begin{definition}
A tensor field $L$ of type $(1,1)$, i.e., an operator field, is called a
Nijenhuis operator if its Nijenhuis torsion is identically zero, i.e., $N_L
\equiv 0$.
\end{definition}

Let $L$ be a Nijenhuis operator, and let $\chi(t) = \det(t\cdot\mathrm{Id}-L)=
t^n + \sigma_1 t^{n-1} + \ldots + \sigma_n$ be its characteristic polynomial.
In \cite{BMK1}, a theorem is proved on the form of the matrix of the Nijenhuis
operator with functionally independent coefficients of the characteristic
polynomial (invariants):
\begin{equation} \label{formula-L}
L = J^{-1} \tilde{L} J,
\end{equation}
where
$$
\tilde{L} =
\begin{bmatrix}
-\sigma_1 & 1 & 0 & \ldots & 0 \\
-\sigma_2 & 0 & 1 & \ldots & 0 \\
\ldots & \ldots & \ldots & \ldots & \ldots \\
-\sigma_n & 0 & 0 & \ldots & 0
\end{bmatrix},
\qquad J = \Bigl(\frac{\partial \sigma_i(x)}{\partial x^j}\Bigr).
$$
For example, in the two-dimensional case, if the trace $\tr L = x$ is the first
coordinate, and $\det L = f(x,y)$ is a smooth function, then $\sigma_1=-x$,
$\sigma_2=f(x,y)$ and, from formula \eqref{formula-L}, we obtain the general
form of these operators:
\begin{equation} \label{L-2dim}
L =
\begin{bmatrix}
%% f_y & 0 \\
%% -f_x & 1
 -1 & 0 \\
\frac{f_x}{f_y}& \frac1{f_y}
\end{bmatrix}
\begin{bmatrix}
x & 1 \\
-f & 0
\end{bmatrix}
\begin{bmatrix}
-1 & 0 \\
f_x & f_y
\end{bmatrix}
= \begin{bmatrix}
x - f_x & -f_y \\
\frac{-xf_x + f^2_x + f}{f_y} & f_x
\end{bmatrix},
\end{equation}
where $f_x,f_y$ stand for the corresponding partial derivatives of the function
$f$. It is clear that the existence of a smooth two-dimensional Nijenhuis
operator in the situation under consideration is equivalent to the smoothness
of the fraction $\frac{xf_x - f^2_x - f}{f_y}$. At the points at which $f_y \ne
0$, this fraction is obviously smooth. A more interesting case occurs when $f_y
= 0$ at some point, but $\frac{xf_x - f^2_x - f}{f_y}$ is smooth in its
neighborhood. The Morse and cubic singularities of the determinant (restricted
to the level line of the trace) in the two-dimensional case were investigated
in \cite{Akpan}.

The following assertion gives an analog of formula \eqref{L-2dim} for the
multidimensional case.

\begin{theorem} \label{th1}
Let $L$ be an $n$-dimensional Nijenhuis operator for which  $\sigma_1 = x_1$,
$\sigma_2 = x_2$, \ldots, $\sigma_{n-1} = x_{n-1}$ in some coordinates
$(x_1,x_2,\dots,x_{n-1},y)$, and $\sigma_n= (-1)^n\det L = f(x_1, \ldots,
x_{n-1}, y)$ is an arbitrary smooth function such that $f_y \ne 0$. Then, in
these coordinates, $L$ has the following form\/\textup:
$$
L = \begin{bmatrix}
-x_1 & 1 & 0 & \dots & 0 & 0 \\
-x_2 & 0 & 1 & \dots & 0 & 0 \\
\dots & \dots & \dots & \dots & \dots & \dots \\
-x_{n-2} & 0 & 0 & \dots & 1 & 0 \\
-x_{n-1} + f_{x_1} & f_{x_2} & f_{x_3} & \dots & f_{x_{n-1}} & f_{y} \\
\frac{x_1 f_{x_1} + x_2 f_{x_2} + \dots + x_{n-1}f_{x_{n-1}} - f_{x_1}
f_{x_{n-1}} - f}{f_y} & -\frac{f_{x_1} + f_{x_2}f_{x_{n-1}}}{f_y} &
-\frac{f_{x_2} + f_{x_3}f_{x_{n-1}}}{f_y} & \dots & -\frac{f_{x_{n-2}} +
f^2_{x_{n-1}}}{f_y} & -f_{x_{n-1}}
\end{bmatrix}.
$$
\end{theorem}

\begin{remark}
In Theorem \ref{th1}, the first $n-2$ rows of the matrix look like a
differential nondegenerate Nijenhuis operator (see \cite{BMK1}), and the
differences are only in the last two lines.
\end{remark}

{\bf Proof.} It is necessary to prove that equality~\eqref{formula-L} is
satisfied, where
$$
J =
\begin{bmatrix}
1 & 0 & \dots & 0 & 0 \\
0 & 1 & \dots & 0 & 0 \\
\dots & \dots & \dots & \dots & \dots \\
0 & 0 & \dots & 1 & 0 \\
f_{x_1} & f_{x_2} & \dots & f_{x_{n-1}} & f_y
\end{bmatrix},
\qquad \tilde{L} =
\begin{bmatrix}
-x_1 & 1 & 0 & \dots & 0 \\
-x_2 & 0 & 1 & \dots & 0 \\
\dots & \dots & \dots & \dots & \dots \\
-x_{n-1} & 0 & 0 & \ldots & 1 \\
-f   & 0 & 0 & \dots & 0
\end{bmatrix}.
$$
To this end, it suffices to show that $JL = \tilde{L}J$. Calculating, we see
that the left- and right-hand sides of this equality coincide:
$$
JL =
\begin{bmatrix}
-x_1 & 1 & 0 & \dots & 0 & 0 \\
-x_2 & 0 & 1 & \dots & 0 & 0 \\
\dots & \dots & \dots & \dots & \dots & \dots \\
-x_{n-2} & 0 & 0 &\dots & 1 & 0 \\
-x_{n-1} + f_{x_1} & f_{x_2} & f_{x_3} & \dots & f_{x_{n-1}} & f_{y} \\
-f & 0 & 0 & \dots & 0 & 0
\end{bmatrix}
= \tilde{L}J.
$$
This completes the proof of the theorem.

\begin{remark} \label{rem2}
It follows from Theorem 1 that the existence of a smooth operator field of the
indicated form is equivalent to the smoothness of the following fractions:
\begin{equation*}
%\begin{cases}
\frac{\sum_{i=1}^{n-1} x_if_{x_i} - f_{x_1} f_{x_{n-1}} - f}{f_y}, \qquad
\frac{f_{x_{j-1}} + f_{x_j}f_{x_{n-1}}}{f_y}, \quad j = 2, \ldots, n-1.
%\end{cases}
\end{equation*}
\end{remark}

In Theorem \ref{th1}, it is actually assumed only that the first $n-1$
coefficients of the characteristic polynomial are functionally independent.
Then they can be taken as the first $n-1$ coordinates. If, in addition, the
Jacobi matrix
$\frac{\partial(\sigma_1,\dots,\sigma_n)}{\partial(x_1,\dots,x_{n-1},y)}$ is
degenerate, then this Nijenhuis operator can be called almost differentially
nondegenerate.

Let us now study the case in which such a singularity of the Jacobi matrix is
given by a nondegenerate singularity of the last coefficient (the determinant)
$\sigma_n = \pm \det L$ with respect to the remaining coordinate $y$.

\begin{theorem} \label{th2}
Let $L$ be an $n$-dimensional Nijenhuis operator, where $n>2$, for which
$\sigma_1 = x_1$, $\sigma_2 = x_2$, \ldots, $\sigma_{n-1} = x_{n-1}$ in some
coordinates $(x_1,x_2,\dots,x_{n-1},y)$, and let $\sigma_n= (-1)^n\det L =
f(x_1, \ldots, x_{n-1}, y)$ be a smooth function with a Morse singularity with
respect to the variable $y$, $f(0, \ldots, 0) = 0$. Then there is a regular change of coordinates
$y\to y(x_1,\dots,x_{n-1},y)$ preserving the remaining coordinates
$x_1,\dots,x_{n-1}$, after which $f=\pm y^2$, and the corresponding Nijenhuis
operators have the form
$$
L =
\begin{bmatrix}
-x_1 & 1 & 0 & \dots & 0 & 0 \\
-x_2 & 0 & 1 & \dots & 0 & 0 \\
\dots & \dots & \dots & \dots & \dots & \dots \\
-x_{n-2} & 0 & 0 & \dots & 1 & 0 \\
-x_{n-1} & 0 & 0 & \dots & 0 & \pm 2y \\
-\frac{y}{2} & 0 & 0 & \dots & 0 & 0
\end{bmatrix}.
$$
\end{theorem}

{\bf Proof.} Let us use the parametric Morse lemma and reduce the function $f$
to the canonical form $f = \pm y^2 + R(x_1,\dots,x_{n-1})$. (Below, by $R_i$,
we denote the partial derivatives $\frac{\partial R}{\partial x_i}$.) A
necessary and sufficient condition for the existence of the Nijenhuis operator
in this case is the smoothness of the following fractions (see Remark
\ref{rem2}):
\begin{multline*}
    \frac{x_1R_1 + x_2 R_2 + \ldots + x_{n-1}R_{n-1} - R_1 R_{n-1} - R \mp y^2}{2y}
= \frac{x_1R_1 + x_2 R_2 + \ldots + x_{n-1}R_{n-1} - R_1 R_{n-1} - R}{2y}  \\ \mp \frac{y}{2},
\end{multline*}
$$
\begin{gathered}
\frac{R_1 + R_2 R_{n-1}}{2y}, \qquad \frac{R_2 + R_3 R_{n-1}}{2y}, \qquad
\dots, \qquad \frac{R_{n-2} + R^2_{n-1}}{2y}.
\end{gathered}
$$

Since, in all cases, we obtain a fraction of the form
$\frac{F(x_1,\dots,x_{n-1})}{y}$, it follows that its smoothness is equivalent
to the condition $F \equiv 0$, i.e., we obtain a system of partial differential
equations with respect to the function $R = R(x_1, \dots, x_{n-1})$:
\begin{equation} \label{pde1}
\begin{cases}
x_1R_1 + x_2 R_2 + \ldots 
+ x_{n-1}R_{n-1} - R_1 R_{n-1} -  R = 0\\
R_1 + R_2 R_{n-1} = 0 \\
R_2 + R_3 R_{n-1} = 0 \\
\ldots \\
R_{n-2} + R^2_{n-1} = 0 \\
\end{cases}
\end{equation}

Note that the equations of this system (starting from the second one) can be
rewritten in the form of the following relations:
\begin{equation} \label{pde2}
R_{n-i} = (-1)^{i-1} R^i_{n-1}, \qquad i = 2, \ldots, n-1.
\end{equation}
Differentiating these relations with respect to $x_1$, we obtain
\begin{equation} \label{pde3}
R_{1,n-i} = i (-1)^{i-1} R^{i-1}_{n-1} R_{1,n-1}, \qquad i = 2, \ldots, n-1,
\end{equation}
where $R_{i,j}$ stand for the corresponding second partial derivatives of the
function $R$. Let us now differentiate the first equation of system
\eqref{pde1} with respect to $x_1$. We obtain the relation
$$
x_1 R_{1,1} + x_2 R_{1,2} + \ldots + x_{n-1}R_{1,n-1} - R_{1,1} R_{n-1} - R_1
R_{1, n-1} = 0
$$
and substitute the expressions for $R_{1,k}$ from formulas \eqref{pde3} into
this relation:
\begin{equation} \label{pde4}
R_{1, n-1}\bigr(x_1 (n-1)(-1)^{n-2} R_{n-1}^{n-2} + \ldots + x_{n-2} 2 (-1)
R_{n-1} + x_{n-1}\bigl) - (-1)^n n R^{n-1}_{n-1} R_{1,n-1} = 0.
\end{equation}

The left-hand side of relation \eqref{pde4} can be factorized, and one of the
factors is $R_{1, n-1}$. Let us study the case in which $R_{1, n-1} = 0$. It
follows from relation \eqref{pde3} that, in this case, $R_{1,1} = R_{1,2} =
\ldots = R_{1,n-1} = 0$, i.e., $R_1 \equiv\mathrm{const}$. Further, from the
fact that $R_1 \equiv\mathrm{const}$, taking into account relations
\eqref{pde2}, we see that all $R_i$ are constant. This means that $R = C_1 x_1
+ \ldots + C_{n-1}x_{n-1}$, where $C_i$ are some constants. Substitute this
expression for $R$ into the first equation of system \eqref{pde1}:
$$
x_1R_1 + x_2 R_2 + \ldots + x_{n-1}R_{n-1} - R_1 R_{n-1} -  R = 0 \quad
\Rightarrow \quad R_1 R_{n-1} = C_1 C_{n-1} = 0.
$$
Each of the conditions $C_1 = 0$ and $C_{n-1} = 0$ implies that the other
constants $C_i$ vanish, since $R_{k-1} + R_k R_{n-1} = 0$, $k = 2, \ldots,
n-1$, in system \eqref{pde1}. Thus, the first solution of the system is $R =
0$.

Let us now show that there are no other smooth solutions. To this end, consider
the second factor in relation \eqref{pde4}:
$$
(-1)^{n-1} n R^{n-1}_{n-1} + x_1(n-1)(-1)^{n-2} R^{n-2}_{n-1} + \ldots +
x_{n-2} 2 (-1) R_{n-1} + x_{n-1} = 0.
$$
Let us rewrite this equation by making the reverse substitution according to
formulas \eqref{pde2}:
\begin{equation} \label{pde5}
nR_{1} + (n-1) x_1 R_2 + \ldots + 3x_{n-3} R_{n-2} + 2x_{n-2} R_{n-1} - x_{n-1}
= 0.
\end{equation}
It is clear from \eqref{pde5} that $R_1(0) = 0$, and then it follows from
system \eqref{pde1} that $R_j(0) = 0$ for $1 < j < n$. Differentiating
\eqref{pde5} with respect to $x_{n-1}$, we obtain $nR_{1,n-1}(0) - 1 = 0$, i.e.,
$R_{1,n-1}(0) = \frac{1}{n}$. On the other hand, differentiating the relation
$R_1 = (-1)^{n-2} R^{n-1}_{n-1}$ with respect to $x_{n-1}$ at zero, we obtain
$R_{1,n-1}(0)=(-1)^{n-2} (n-1)R^{n-2}_{n-1}(0) R_{n-1, n-1}(0) = 0$. We arrive
at a contradiction and hence such a smooth function does not exist.

This completes the proof of Theorem \ref{th2}.

\begin{remark}
The case $n = 2$ is special, and there arises an additional solution $R =
\frac{x^2}{4}$, since, in this case, there are no additional conditions on the
function $R(x)$ except for $xR_x - R^2_x - R = 0$.
\end{remark}

\section*{Acknowledgments}
The author thanks A.\, A.~Oshemkov, E.\, A.~Kudryavtseva, A.\, Yu.~Konyaev,
V.\, A.~Kibkalo, V.\, N.~Zavyalov, and all participants of the seminar
``Algebra and Geometry of Integrable Systems'' for valuable advice.

\end{document}